\documentclass[12pt]{article}
\usepackage{amsmath}
\usepackage{amssymb}
\usepackage{amsbsy}

\usepackage{amsfonts}

\newtheorem{lem}{Lemma}[section]%
\newtheorem{thm}[lem]{Theorem}%
\newtheorem{prop}[lem]{Proposition}%

\def\a{\alpha}  \def\g{\gamma}  
 \def\s{\sigma}

\def\G{\Gamma}

 \def\lg{\langle} \def\rg{\rangle}
\def\nd{\mathrel{\bigm|\kern-.7em/}}

\def\f{\noindent}

\def\Aut{\hbox{\rm Aut\,}}

\def\mod{\hbox{\rm mod }}

\def\mz{{\mathbb Z}}
\def\C{\hbox{\rm C}}

\def\demo{\noindent{\bf Proof}\hskip10pt}
\def\case{\noindent{\bf Case}\hskip5pt}
\def\qed{\hskip10pt $\Box$}
\def\qqed{\hskip10pt $\Box$\vspace{3mm}}

\begin{document}
\title{Absolutely split metacyclic groups and weak metacirculants}
\author
{Li Cui, Jin-Xin Zhou   \\
{\small Department of Mathematics, Beijing Jiaotong University}\\ {\small Beijin 100044, P.R. China}\\
{\small{Email}:\ licui@bjtu.edu.cn, jxzhou@bjtu.edu.cn}}

\date{}
\maketitle

\begin{abstract}
Let $m,n,r$ be positive integers, and let $G=\lg a\rg: \lg b\rg\cong\mz_n:\mz_m$ be a split metacyclic group such that $b^{-1}ab=a^r$. We say that $G$ is {\em absolutely split with respect to $\lg a\rg$} provided that for any $x\in G$, if $\lg x\rg\cap\lg a\rg=1$, then there exists $y\in G$ such that $x\in\lg y\rg$ and $G=\lg a\rg:\lg y\rg$. In this paper, we give a sufficient and necessary condition for the group $G$ being absolutely split. This generalizes a result of Sanming Zhou and the second author in [arXiv: 1611.06264v1]. We also use this result to investigate the relationship between metacirculants and weak metacirculants.

Metacirculants were introduced by Alspach and Parsons in $1982$ and have been a rich source of various topics since then. As a generalization of this classes of graphs, Maru\v si\v c and \v Sparl in 2008 posed the so called weak metacirculants. A graph is called a {\em weak metacirculant} if it has a vertex-transitive metacyclic automorphism group. In this paper, it is proved that a weak metacirculant of $2$-power order is a metacirculant if and only if it has a vertex-transitive split metacyclic automorphism group. This provides a partial answer to an open question in the literature.

\bigskip
\noindent{\bf Key words:} 2-group, absolutely split metacyclic group, metacirculant.\\
\noindent{\bf 2000 Mathematics subject classification:} 20B25, 05C25.
\end{abstract}
\section{Introduction}
A group $G$ is called {\em metacyclic} if it contains a cyclic normal subgroup $N$ such that $G/N$ is cyclic. In other words, a metacyclic group $G$ is an extension of a cyclic group  $N\cong C_n$ by a cyclic group $G/N\cong C_m$, written $G\cong C_n.C_m$. If this extension is split, namely $G \cong C_n: C_m$, then $G$ is called a {\em split metacyclic group}. Metacyclic groups form a basic and well-studied family of groups. Certain classes of metacyclic groups have been given much attention, see, for example, \cite{Hempel,King,Liedahl,Newman-Xu,Schenkman,Sim,Xu-Zhang}.

In this paper, we shall be concerned with the split metacyclic groups which are closed related to the metacirculants. Let $m\geq 1$ and $n\geq 2$ be integers. A graph $\G = (V(\G), E(\G))$ of order $mn$ is called \cite{MS} an {\em $(m, n)$-metacirculant graph} (in short {\em $(m, n)$-metacirculant}) if it has two automorphisms $\s,\tau$ such that
\begin{enumerate}
  \item [{\rm (1)}]\ $\lg \s\rg$ is semiregular and has $m$ orbits on $V(\G)$,
  \item [{\rm (2)}]\ $\tau$ normalizes $\lg\s\rg$ and cyclically permutes the $m$ orbits of $\lg\s\rg$,
  \item [{\rm (3)}]\ $\tau$ has a cycle of size $m$ in its cycle decomposition.
\end{enumerate}
A graph is called a {\em metacirculant} if it is an $(m, n)$-metacirculant for some $m$ and $n$.
Metacirculant graphs were introduced by Alspach and Parsons \cite{AP} in 1982, and have attracted a lot of attention. It follows from the definition above that a metacirculant $\G$ has an autormorphism group $\lg \s, \tau\rg$ which is metacyclic and transitive on $V(\G)$.

As a generalization of metacirculants, Maru\v si\v c and \v Sparl~\cite{MS} posed the so called weak metacirculants. A graph is called a {\em weak metacirculant} if it has a vertex-transitive metacyclic automorphism group. In \cite{LWS} Li et al. divided the metacirculants into the following two subclasses: A weak metacirculant which has a vertex-transitive split metacyclic automorphism group is called {\em split weak metacirculant}.
Otherwise, a weak metacirculant $\Gamma$ is called a {\em non-split weak metacirculant}
if its full automorphism group does not contain any split metacyclic subgroup which is vertex-transitive.
In \cite{LWS} Li et al. studied the relationship between metacirculants and weak metacirculants.
Among other results they proved that every metacirculant is a split weak metacirculant (see \cite[Lemma~2.2]{LWS}), but it was unknown whether the converse of this statement is true. In \cite[Question~A]{Zhou} Sanming Zhou and the second author asked the following question:

\medskip
\f{\bf Question~A}\, Is it true that any split weak metacirculant is a metacirculant?
\medskip

In the study of the relationship between metacirculants and weak metacirculants, the so called absolutely split metacyclic groups (defined below) play an important role. Let $m,n,r$ be positive integers, and let $G=\lg a\rg: \lg b\rg\cong\mz_n:\mz_m$ be a split metacyclic group such that $b^{-1}ab=a^r$. We say that $G$ is {\em absolutely split with respect to $\lg a\rg$} provided that for any $x\in G$, if $\lg x\rg\cap\lg a\rg=1$, then there exists $y\in G$ such that $x\in\lg y\rg$ and $G=\lg a\rg:\lg y\rg$. We say that a split metacyclic group is {\em absolutely split} if it is absolutely split with respect to its some normal cyclic subgroup.

Clearly, if a connected weak metacirculant has a vertex-transitive absolutely split metacyclic automorphism group, then it is also a metacirculant.
Actually, by proving that the group $G\cong\mz_{p^n}:\mz_{p^m}$ with $p$ an odd prime and $n\geq m\geq 1$ is absolutely split, Zhou and the second author in \cite[Theorem~1.1]{Zhou} proved that a connected weak metacirculant with order an odd prime power is a metacirculant if and only if it is a split weak metacirculant.

In this paper, we give a sufficient and necessary condition for a split metacyclic group being absolutely split, and this is then used to prove that a connected weak metacirculant with order a $2$-power is a metacirculant if and only if it is a split weak metacirculant. This together with \cite[Theorem~1.1]{Zhou} shows that the answer to Question~A is positive when the graph under consideration is of order a prime power.

However, in general the answer to Question~A is negative, and an infinite family of split weak metacirculants which are not metacirculants will be constructed in our subsequent paper~\cite{Cui-Zhou}.


\section{Definitions and notations}
For a positive integer $n$, we denote by $C_{n}$ the cyclic group of order $n$,
by $\mathbb{Z}_{n}$ the ring of integers modulo $n$,
and by $\mathbb{Z}_{n}^{\ast}$ the multiplicative group of $\mathbb{Z}_{n}$ consisting of numbers coprime to $n$.

Let $G$ be a finite group.
The full automorphism group, the center, the derived group and the Frattini subgroup of $G$
will be denoted by ${\Aut(G)}$, $Z(G)$, $G'$ and $\Phi(G)$, respectively.
For $x\in G$, denote by $o(x)$ the order of $x$.
For $a_{1}, a_{2},\ldots, a_{n}\in G$,
$n\ge 2$,
the commutator $[a_{1}, a_{2},\ldots, a_{n}]$ of $a_{1}, a_{2},\ldots, a_{n}$ was recursively defined as follow:
if $n=2$, then $[a_{1},a_{2}]=a_{1}^{-1}a_{2}^{-1}a_{1}a_{2}$; and
if $n>2$, then $[a_{1}, a_{2},\ldots, a_{n}]=[[a_{1}, a_{2},\ldots, a_{n-1}],a_{n}].$

Let $p$ be a prime and let $G$ be a $p$-group of exponent $p^{e}$.
For any $0\le s\le e$, let
$$\Omega_{s}(G)=\lg g\in G \mid g^{p^{s}}=1\rg.$$

For a finite, simple and undirected graph $\Gamma$,
we use $V(\Gamma)$, $E(\Gamma)$, ${\rm Aut(\Gamma)}$ to denote its vertex set, edge set and full automorphism group, respectively.
Let $G\leq {\rm Aut(\Gamma)}$, $v\in V(\Gamma)$.
Denote by $G_{v}$ the stabilizer of $v$ in $G$,
that is, the subgroup of $G$ fixing the point $v$.
We say that $G$ is {\em semiregular} on $V(\Gamma)$ if $G_{v}=1$ for every $v\in V(\Gamma)$
and {\em regular} if $G$ is transitive and semiregular.

\section{A technical lemma}

A group $G$ is said to be {\em metabelian group} if $G''=1$,
that is the derived group $G'$ is abelian.
For a metabelian group $G$,
if $x\in G'$ and $x_{1}, x_{2}, \ldots , x_{n}\in \{ a, b\}$
and if $\sigma$ is a permutation on the set $\{ 1,2, \ldots, n\}$,
then $[x, x_{1}, x_{2}, \ldots, x_{n}]=[x, x_{1^{\sigma}}, x_{2^{\sigma}}, \ldots , x_{n^{\sigma}}]$.
So for brevity of writing we make the following convention:
$$[ia,jb]=[a,b,\underbrace{a,\ldots,a}_{i-1},\underbrace{b,\ldots,b}_{j-1}].$$

The following result is due to Xu.

\begin{prop}{\rm\cite[Lemma~3]{Xu}}\label{2}
Let $G$ be a metabelian group and let $\ell\geq2$ be an integer. Then for any $x,y\in G$, we have
$$(xy^{-1})^{\ell}=x^{\ell}(\prod_{i+j\le \ell}[ix,jy]^{\C_{\ell}^{i+j}})y^{-\ell}.$$
$($Here for any integers $N\geq l\geq 0$, we denote by $\C_N^l$ the binomial coefficient, that is, $\C_N^l=\frac{N!}{l!(N-l)!}$.$)$
\end{prop}

Using Proposition~\ref{2}, we have the following lemma which will be frequently used in the following sections.

\begin{lem}\label{3}
Let $G=\lg a\rg.\lg b\rg\cong C_{n}.C_{m}$.
Then for any integers $k\ge 2,$ and $i,j\geq 0$, we have
$$(b^{i}a^{j})^{k}=b^{ik}a^{jk}\prod_{2\le s+1\le k}[a,sb^{i}]^{j\C_{k}^{s+1}},$$
and moreover, if $a^b=a^r$ for some $r\in\mz_n$, then
$$(b^{i}a^{j})^{k}
=b^{ik}a^{j(k+\sum_{2\le s+1\le k}(r^{i}-1)^{s}\C_{k}^{s+1})}.$$
\end{lem}

\demo Observing that $\lg a \rg\unlhd G$ and $G/\lg a\rg\cong\mz_m$, one has $G'\leq \lg a\rg$.
So for every $g\in G'$ and $h\in\lg a\rg$, we have $[g,h]=[h,g]=1.$ We first prove two claims.\medskip

\f{\bf Claim~1}\ {\em For any $g_1\in G, g_2\in\lg a\rg$ and for any integer $i$, we have $[g_{1},g_{2}^{i}]=[g_{1},g_{2}]^{i}$ and $[g_{2}^{i},g_{1}]=[g_{2},g_{1}]^{i}$.} \medskip

We shall only prove the first formula. The second one can be proved in a similar way.  Assume that $i\geq 1$. If $i=1$, our claim is clearly true. Assume that $i>1$.  Then $$[g_{1}, g_{2}^{i}]=[g_{1},g_{2}^{i-1}g_{2}]=[g_{1},g_{2}][g_{1},g_{2}^{i-1}]^{g_{2}}=[g_{1},g_{2}][g_{1},g_{2}^{i-1}].$$
By induction on $i$, we have $$[g_{1}, g_{2}^{i}]=[g_{1},g_{2}]^{i}.$$
Noting that $g_{2}^{-i}=g_{2}^{o(g_{2})-i}=g_{2}^{o(a)-i}$, it follows that
$$[g_{1},g_{2}^{-i}]=[g_{1},g_{2}^{o(a)-i}]=[g_{1},g_{2}]^{o(a)-i}=[g_{1},g_{2}]^{-i}.$$
Hence $[g_{1},g_{2}^{i}]=[g_{1},g_{2}]^{i}$ holds for any integer $i$.
\medskip

\f{\bf Claim~2}\ {\em For any positive integer $s$ and for any integers $i,j$, we have $[sb^{i},a^{-j}]=[a,sb^{i}]^{j}.$}\medskip



Repeatedly using Claim~1, we have
\begin{align*}
[sb^{i},a^{-j}]=&[[b^{i},a^{-j}],(s-1)b^{i}]\\
=&[[a,b^{i}]^{j},(s-1)b^{i}]\\
=&[[[a,b^{i}]^{j},b^{i}],(s-2)b^{i}]\\
=&[[[a,b^{i}],b^{i}]^{j},(s-2)b^{i}]\\
=&[[a,2b^{i}]^{j},(s-2)b^{i}]\\
\ldots&\\
=&[a,sb^{i}]^{j}.
\end{align*}

Now we are ready to complete the proof of our lemma. For the first part, by Proposition~\ref{2}, for any integers $k\ge 2,$ and $i,j\geq 0$,  we have
\begin{align*}
(b^{i}a^{j})^{k}=& b^{ik}\prod_{2\le s+t\le k}[sb^{i},ta^{-j}]^{\C_{k}^{s+t}}a^{jk}\\
=& b^{ik}a^{jk}\prod_{2\le s+t\le k}[sb^{i},ta^{-j}]^{\C_{k}^{s+t}}.\\
\end{align*}
Clearly, $G'\leq\lg a\rg$, so if $t\ge 2$, we have $[sb^{i},ta^{-j}]=[sb^{i},(t-1)a^{-j},a^{-j}]=1.$
Combining this formula with Claim~2, we have
\begin{align*}
b^{ik}a^{jk}\prod_{2\le s+t\le k}[sb^{i},ta^{-j}]^{\C_{k}^{s+t}}=& b^{ik}a^{jk}\prod_{2\le s+1\le k}[sb^{i},a^{-j}]^{\C_{k}^{s+1}}\\
=& b^{ik}a^{jk}\prod_{2\le s+1\le k}[a,sb^{i}]^{j\C_{k}^{s+1}}.
\end{align*}
This establishes the first formula of our lemma.

Now let $a^{b}=a^{r}$ for some $r\in\mz_n$. Then
$$[a,b^{i}]=a^{-1}b^{-i}ab^{i}=a^{-1}a^{b^{i}}=a^{-1}a^{r^{i}}=a^{r^{i}-1}.$$
If $s\ge 2$, then by Claim~1, we have
$$[a, sb^{i}]=[[a,b^{i}],(s-1)b^{i}]=[a^{r^{i}-1},(s-1)b^{i}]=[a,(s-1)b^{i}]^{r^{i}-1}.$$
By induction on $s$, we get the following formula $$[a,sb^{i}]=[a,b^{i}]^{(r^{i}-1)^{s-1}}=a^{(r^{i}-1)^{s}}.$$
Therefore,
\begin{align*}
(b^{i}a^{j})^{k}
=& b^{ik}a^{jk}\prod_{2\le s+1\le k}[a,sb^{i}]^{j\C_{k}^{s+1}}\\
=& b^{ik}a^{jk}\prod_{2\le s+1\le k}(a^{(r^{i}-1)^{s}})^{j\C_{k}^{s+1}}\\
=& b^{ik}a^{j[k+\sum_{2\le s+1\le k}(r^{i}-1)^{s}\C_{k}^{s+1}]}.
\end{align*}
\hfill\qed

\section{Absolutely split metacyclic groups}

In this section, we shall give a characterization of absolutely split metacyclic groups. Throughout this section, we shall make the following assumptions:\medskip

\f{\bf Assumption.}\
\begin{itemize}
\item $n, m$:\ two positive integers;
\item $G=\lg a\rg:\lg b\rg\cong C_{n}:C_{m}$, where $b^{-1}ab=a^r$ for some $1\neq r\in\mz_n^*$.
\end{itemize}

Note that every element of $G$ can be written as the form $b^{i}a^{j}$ for some $i\in\mz_m, j\in\mz_n$. A pair $(i, j)\in\mz_m\times\mz_n$ is said to be {\em admissible} with respect to $\lg a\rg$ if $\lg a\rg\cap\lg b^{i}a^{j}\rg=1$.

\begin{lem}\label{lem-admissible}
A pair $(i, j)\in\mz_m\times\mz_n$ is admissible if and only if $o(b^{i}a^{j})=o(b^{i})$.
\end{lem}

\demo Let $x=b^{i}a^{j}$. By Lemma~\ref{3}, for every $k\geq 2$, we have
$$x^k=(b^{i}a^{j})^{k}=b^{ik}a^{jk}\prod_{2\le s+1\le k}[a,sb^{i}]^{j\C_{k}^{s+1}}.$$
Since $G/\lg a\rg$ is cyclic, one has $G'\leq\lg a\rg$, and hence $$a^{jk}\prod_{2\le s+1\le k}[a,sb^{i}]^{j\C_{k}^{s+1}}\in \lg a\rg.$$
It follows that
\begin{equation}\label{eq-1}
(b^{ik})^{-1}x^k=(b^{ik})^{-1}(b^{i}a^{j})^{k}\in \lg a\rg.
\end{equation}

Suppose first that $(i, j)\in\mz_m\times\mz_n$ is admissible. Then $\lg a\rg\cap\lg x\rg=1$. If $x=1$, then $b^{i}=a^{-j}\in\lg a\rg\cap\lg b\rg=1$, and so $b^i=1$, as required. If $b^i=1$, then $x=a^j\in\lg a\rg\cap\lg x\rg=1$, and so $x=1$, as required. So we may assume that $x$ is non-trivial. Applying Eq.~(\ref{eq-1}), we obtain that $x$ and $b^i$ have the same order due to $\lg x\rg\cap\lg a\rg=\lg b\rg\cap \lg a\rg=1$.

Suppose now that $o(x)=o(b^i)=t$. We may assume that $t>1$. Take $x^k\in\lg a\rg\cap\lg x\rg$. If $k=1$, then $x\in\lg a\rg$, and so $b^i=xa^{-j}\in\lg a\rg\cap\lg b\rg=1$. Consequently, we have $b^i=1$, and so $x=1$, as required. Suppose that $k>1$.  Again by Eq.~(\ref{eq-1}), we have $b^{ik}\in \lg a\rg\cap\lg b\rg=1$, and hence $b^{ik}=1$. It follows that $t\ |\ k$, and so $x^k=1$. Therefore, we have $\lg a\rg\cap\lg x\rg=1$, and hence $(i, j)\in\mz_m\times\mz_n$ is admissible.\hfill\qed

\begin{lem}\label{lem-admissible-2}
A pair $(i, j)\in\mz_m\times\mz_n$ is admissible if and only if the following equation
\begin{equation}\label{eq-admissible-2-0}
j[r^{i(k-1)}+r^{i(k-2)}+\cdots+r^i+1]\equiv 0\ (\mod n),
\end{equation}
holds, where $k$ is the order of $b^{i}$.
\end{lem}

\demo If $k=1$, the lemma is clearly true. In what follows, we assume that $k>1$. By Lemma~\ref{3}, we have
\begin{equation}\label{eq-admissible-2}
(b^{i}a^{j})^{k}=b^{ik}a^{jk+\sum_{2\le s+1\le k}\C_{k}^{s+1}(r^{i}-1)^{s}j}=a^{jk+\sum_{2\le s+1\le k}\C_{k}^{s+1}(r^{i}-1)^{s}j}.
\end{equation}

By Lemma~\ref{lem-admissible}, the pair $(i, j)\in\mz_m\times\mz_n$ is admissible if and only if $o(b^{i})=o(b^{i}a^{j})=k.$
By Eq.~(\ref{eq-admissible-2}), if $o(b^{i}a^{j})=k$, then we have
\begin{equation}\label{eq-admissible-2-1}
a^{j[k+\sum_{2\le s+1\le k}\C_{k}^{s+1}(r^{i}-1)^{s}]}=1.
\end{equation}
Conversely, if Eq.~(\ref{eq-admissible-2-1}) holds, then Eq.~(\ref{eq-admissible-2}) implies that $(b^{i}a^{j})^{k}=1$, and so $o(b^{i}a^{j})\ |\ k$. Let $\ell=o(b^{i}a^{j})$. Since $k>1$, one has $\ell>1$, and again by
Lemma~\ref{3}, we have
\[
1=(b^{i}a^{j})^{\ell}=b^{i\ell}a^{j[\ell+\sum_{2\le s+1\le \ell}\C_{\ell}^{s+1}(r^{i}-1)^{s}]}.
\]
It then follows that $b^{i\ell}=1$. Hence $o(b^i)=k$ divides $\ell$, and consequently, we have $o(b^{i})=o(b^{i}a^{j})$.
Now we conclude that the pair $(i, j)\in\mz_m\times\mz_n$ is admissible if and only if Eq.~(\ref{eq-admissible-2-1}) holds.

To complete the proof, it suffices to show that Eq.~(\ref{eq-admissible-2-1}) is equivalent to Eq.~(\ref{eq-admissible-2-0}).
Clearly, Eq.~(\ref{eq-admissible-2-1}) is equivalent to the following equation
\[j[k+\sum_{2\le s+1\le k}\C_{k}^{s+1}(r^{i}-1)^{s}]\equiv 0\ (\mod n).\]
If $r^i=1$, then the above equation is just Eq.~(\ref{eq-admissible-2-0}).
Suppose now that $r^i\neq 1$. Multiplying the above equation by $r^{i}-1$, we have
\[j[k(r^{i}-1)+\sum_{2\le s+1\le k}\C_{k}^{s+1}(r^{i}-1)^{s+1}]\equiv 0\ (\mod n(r^{i}-1)),\]
namely,
\[j[\sum_{0\le s+1\le k}\C_{k}^{s+1}(r^{i}-1)^{s+1}-1]\equiv 0\ (\mod n(r^{i}-1)).\]
Then we have \[j(r^{ik}-1)\equiv 0\ (\mod n(r^{i}-1)).\]
Since $r^i\neq 1$, Equation~(\ref{eq-admissible-2-0}) can be obtained by dividing the above equation by $r^{i}-1$. Therefore, Eq.~(\ref{eq-admissible-2-1}) is equivalent to Eq.~(\ref{eq-admissible-2-0}).
\hfill\qed

\begin{thm}\label{thm-absolute-1}
The metacyclic group $G$ is absolutely split with respect to $\lg a\rg$ if and only if for any admissible pair $(i, j)\in\mz_m\times\mz_n$, there exists an admissible pair $(1, t)\in\mz_m^*\times\mz_n$ such that
\begin{equation}\label{eq-absolute-1}
j\equiv \frac{t(r^{i}-1)}{r-1}\ (\mod n).
\end{equation}
\end{thm}

\f\demo By the definition, $G$ is absolutely split with respect to $\lg a\rg$ if and only if for any admissible pair $(i, j)\in\mz_m\times\mz_n$, there exists an admissible pair $(s, t')\in\mz_m^*\times\mz_n$ such that $b^ia^j=(b^sa^{t'})^\ell$ for some $1\leq \ell\leq m-1$.

Consider the necessity of this assertion. Let $(i, j)\in\mz_m\times\mz_n$ be an admissible pair. Suppose that there exists an admissible pair $(s, t')\in\mz_m^*\times\mz_n$ such that $b^ia^j=(b^sa^{t'})^\ell$ for some $1\leq \ell\leq m-1$. Since $s\in\mz_{m}^*$, there exists $s^{-1}\in \mz_m^*$ such that $ss^{-1}\equiv 1~ (\mod m)$.
By Lemma~\ref{3}, we have $(b^sa^{t'})^{s^{-1}}=ba^{t}$ for some $t\in \mz_n$, and then $\lg ba^t\rg=\lg (b^{s}a^{t'})^{s^{-1}}\rg=\lg b^{s}a^{t'}\rg$. This implies that $(s, t')$ is admissible if and only if $(1, t)$ is admissible. Moreover,
\[b^ia^j=(b^sa^{t'})^\ell=(b^sa^{t'})^{s^{-1}s\ell}=[(b^sa^{t'})^{s^{-1}}]^{s\ell}=(ba^{t})^{s\ell}.\]
By Lemma~\ref{3}, we have $(ba^{t})^{s\ell}=b^{s\ell}a^{t''}$ for some $t''\in \mz_n$.
Thus, $b^ia^j=b^{s\ell}a^{t''}$ and hence $s\ell\equiv i\ (\mod m)$.

The argument in the above paragraph shows that $G$ is absolutely split with respect to $\lg a\rg$ if and only if for any admissible pair $(i, j)\in\mz_m\times\mz_n$, there exists an admissible pair $(1, t)\in\mz_m^*\times\mz_n$ such that $b^ia^j=(ba^{t})^i$.
To finish the proof, it suffices to show the necessity of this assertion is equivalent to Eq.~(\ref{eq-absolute-1}) holds. This is clearly true for the case where $i=1$.


If $i>1$, then by Lemma~\ref{3}, we have
\[(ba^{t})^{i}=b^{i}a^{t[i+\sum_{2\leq k+1\leq i}\C_{i}^{k+1}(r-1)^{k}]},\]
and it then follows that
\begin{align*}
  & b^ia^j=(ba^t)^i\\
  \Longleftrightarrow\ \ & b^{i}a^{j}=b^{i}a^{t[i+\sum_{2\leq k+1\leq i}\C_{i}^{k+1}(r-1)^{k}]} \\
  \Longleftrightarrow\ \  & a^{j}=a^{t[i+\sum_{2\leq k+1\leq i}\C_{i}^{k+1}(r-1)^{k}]}\\
  \Longleftrightarrow\ \  & j\equiv t[i+\sum_{2\leq k+1\leq i}\C_{i}^{k+1}(r-1)^{k}]\ (\mod n)\\
  \Longleftrightarrow\ \  & j(r-1)\equiv t[i(r-1)+\sum_{2\leq k+1\leq i}\C_{i}^{k+1}(r-1)^{k+1}]\ (\mod n(r-1))\\
  \Longleftrightarrow\ \  & j(r-1)\equiv t[(1+(r-1))^{i}-1]\ (\mod n(r-1)) \\
  \Longleftrightarrow\ \  & j(r-1)\equiv t(r^{i}-1)\ (\mod n(r-1))\\
  \Longleftrightarrow\ \  & j\equiv \frac{t(r^{i}-1)}{r-1}\ (\mod n).
\end{align*}
This completes the proof.
\hfill\qed

\medskip
\f{\bf Remark on Theorem~\ref{thm-absolute-1}}\  From the proof of Theorem~\ref{thm-absolute-1}, one may see that if for an admissible pair $(i, j)\in\mz_m\times\mz_n$, there exists an admissible pair $(1, t)\in\mz_m^*\times\mz_n$ satisfying Eq.~(\ref{eq-absolute-1}), then for any admissible pair $(i', j')\in\mz_m\times\mz_n$ such that $\lg b^ia^j\rg=\lg b^{i'}a^{j'}\rg$, there must exist an admissible pair $(1, t')\in\mz_m^*\times\mz_n$ satisfying Eq.~(\ref{eq-absolute-1}).

\section{Absolutely split metacyclic $p$-groups}

We beginning by proving that if a split metacyclic $p$-group is absolutely split with respect to some normal cyclic subgroup of order $n$, then it is also absolutely split with respect to all normal cyclic subgroups of order $n$ which have a complement.

\begin{thm}\label{th-p-groups}
Let $p$ be prime and let $G$ be a split metacyclic $p$-group. Suppose that $G$ has two pair generators $(x,y), (a,b)$ such that $G=\lg x\rg: \lg y\rg=\lg a\rg:\lg b\rg$. If $G$ is non-abelian, then $\lg x\rg\cong\lg a\rg$.
Furthermore, if $G$ is absolutely split with respect to $\lg a\rg$, then $G$ is also absolutely split with respect to $\lg x\rg$.
\end{thm}

\f\demo Since $G/\lg x\rg$ and $G/\lg a\rg$ is abelian,
one has that $1\neq G'\leq \lg x\rg\cap\lg a\rg$.
Since $G$ is a non-abelian $p$-group,
one has $1\neq\Omega_{1}(G')=\Omega_{1}(\lg x\rg)=\Omega_{1}(\lg a\rg).$

Since $\lg y\rg\cap\lg x\rg=1$
and $1\neq\Omega_{1}(\lg x\rg)=\Omega_{1}(\lg a\rg)$,
we have that $\lg y\rg\cap \lg a\rg=1$.
Assume $y=b^{i}a^{j}$.
By Lemma~$4.1$, we have that $o(y)=o(b^{i})$.
Therefore $o(y)\leq o(b)$.
With a similar argument as above we shall have $o(b)\leq o(y)$.
Consequently, $o(b)=o(y)$, and so
$\lg x\rg\cong\lg a\rg$.

Assume that $G$ is absolutely split with respect to $\lg a\rg$. To show that $G$ is also absolutely split with respect to $\lg x\rg$, we take $g\in G$ such that $\lg g\rg\cap\lg x\rg=1$. Recalling that $1\neq \Omega_{1}(\lg a\rg)=\Omega_{1}(\lg x\rg)$, one has $\lg g\rg\cap\lg a\rg=1$,
and then exists $c\in G$ such that $g\in \lg c\rg\cong\lg b\rg$ and
$G=\lg a\rg:\lg c\rg$. Since $o(b)=o(y)$, one has $\lg y\rg\cong\lg b\rg\cong\lg c\rg$. Again, since $1\neq \Omega_{1}(\lg a\rg)=\Omega_{1}(\lg x\rg)$, one has $\lg x\rg\cap\lg c\rg=1$, completing the proof.
\hfill\qed

\medskip
\f{\bf Remark on Theorem~\ref{th-p-groups}}\  Theorem~\ref{th-p-groups} may be not true when $G$ is not a $p$-group.
For example, let $G=(C_{n}: C_{m})\times C_{\ell}$, where $m,n,\ell$ are three positive integers such that $(n, \ell)=(m, \ell)=1$ and $\ell>1$.
Then $G=C_{n\ell}: C_{m}=C_{n}: C_{m\ell}$, but $C_{n\ell}\ncong C_{n}$.
\medskip

In \cite[Lemma~3.4]{Zhou}, it was proved that the group $G\cong\mz_{p^n}:\mz_{p^m}$ with $p$ an odd prime and $n\geq m\geq 1$ is absolutely split.
In this section, we shall consider the  split metacyclic $2$-groups, and prove the following theorem.

\begin{thm}\label{Th}
Let $G$ be a non-abelian split metacyclic $2$-group. If the center of $G$ is cyclic,
then $G$ is absolutely split.
\end{thm}

This theorem will be proved by the following series of lemmas.

\subsection{Split metacyclic $2$-groups with cyclic centers}

To prove Theorem~\ref{Th}, we need the following result,
which is due to Newman, Xu and Zhang.

\begin{prop}{\rm\cite[Theorem 3.1]{Xu-Zhang}}\label{fenl}
Let $G$ be a metacyclic group of order a $2$-power which has no cyclic maximal subgroups.
Then G has one presentation of the following two kinds$:$
\begin{enumerate}

\item [{\rm (I)}]\ Ordinary metacyclic $2$-groups:
\[G=\lg a,b\ |\ a^{2^{r+s+u}}=1, b^{2^{r+s+t}}=a^{2^{r+s}},a^{b}=a^{1+2^{r}}\rg,\]
where $r,s,t,u$ are non-negative integers with $r\ge 2$ and $u\le r$.

\item [{\rm (II)}]\ Exceptional metacyclic $2$-groups:
\[G=\lg a,b\ |\ a^{2^{r+s+v+t'+u}}=1,b^{2^{r+s+t}}=a^{2^{r+s+v+t'}},a^{b}=a^{-1+2^{r+v}}\rg,\]
where $r,s,v,t,t',u$ are non-negative integers with $r\ge 2, t'\le r,u\le 1,tt'=sv=tv=0$, and if $t'\ge r-1$, then $u=0$.
\end{enumerate}
Groups of different types or of the same type but with different values of parameters are not isomorphic to each other.

Furthermore, a Type I group is split if and only if $stu=0$, and
a Type II group is split if and only if $u=0$.

\end{prop}

\begin{lem}\label{G}
Let $G$ be a split metacyclic $2$-group
which has no cyclic maximal subgroups. If $Z(G)$ is cyclic,
then $G$ has a representation$:$
$$G=\lg a,b\ |\ a^{2^{2r+s+t}}=1, b^{2^{r+s}}=1,a^{b}=a^{\pm 1+2^{r+t}}\rg,$$
where $ r\ge 2, st=0$.
\end{lem}

\f\demo
By Proposition~\ref{fenl}, $G$ is a group of Type I or Type II.\medskip

\case 1\ {\em $G$ is a Type I group.}\medskip

In this case, by Proposition~\ref{fenl}, $G$ has a representation$:$
$$G=\lg a,b\ |\ a^{2^{r+s+u}}=1, b^{2^{r+s+t}}=a^{2^{r+s}},a^{b}=a^{1+2^{r}}\rg,$$
where $r,s,t,u$ are non-negative integers with $r\ge 2$, $u\le r$ and $stu=0$.
It is easy to prove that $Z(G)=\lg a^{2^{s+u}}\rg\lg b^{2^{s+u}}\rg$ (see also \cite[p.27]{Xu-Zhang}).
Since $Z(G)$ is cyclic, one has
$Z(G)=\lg a^{2^{s+u}}\rg$ or $ \lg b^{2^{s+u}}\rg.$

If $Z(G)=\lg a^{2^{s+u}}\rg$, then $\lg b^{2^{s+u}}\rg\le \lg a\rg\cap\lg b\rg=\lg b^{2^{r+s+t}}\rg$.
Due to $u\le r$, we also have $\lg b^{2^{r+s+t}}\rg\le\lg b^{2^{s+u}}\rg$.
It then follows that $\lg b^{2^{s+u}}\rg=\lg b^{2^{r+s+t}}\rg$, and hence $u=r$ and $t=0$.

If $Z(G)=\lg b^{2^{s+u}}\rg$, then with a similar argument as above, we may obtain $\lg a^{2^{r+s}}\rg=\lg a^{2^{s+u}}\rg$.
Then $u=r$, and then
$$G=\lg a,b\ |\ a^{2^{2r+s}}=1,b^{2^{r+s+t}}=a^{2^{r+s}},a^{b}=a^{1+2^{r}}\rg,$$ where $r\ge 2, st=0$.
If $t=0$, then $G=\lg a\rg\rtimes\lg ba^{2^{r-1}-1}\rg$, and then by letting $x=a,y=ba^{2^{r-1}-1}$, we have
$$G=\lg x,y\ |\ x^{2^{2r+s}}=1, y^{2^{r+s}}=1, x^{y}=x^{1+2^{r}}\rg.$$
If $t\neq 0,$ then $s=0,G=\lg b\rg\rtimes\lg ab^{-2^{t}}\rg$, and then by letting $x=b,y=ab^{-2^{t}}$, we have
$$G=\lg x,y\ |\ x^{2^{2r+t}}=1, y^{2^{r}}=1, x^{y}=x^{1+2^{r+t}}\rg.$$
Thus, $G$ always has the desired representation.

\medskip
\case 2\ {\em  Assume $G$ is a Type II group.}\medskip

In this case, again by Proposition~\ref{fenl}, $G$ has the following representation$:$
$$G=\lg a,b\ |\ a^{2^{r+s+v+t'}}=1,b^{2^{r+s+t}}=1,a^{b}=a^{-1+2^{r+v}}\rg,$$
where $r,s,v,t,t'$ are non-negative integers with $r\ge 2, t'\le r, tt'=sv=tv=0$.

From \cite[p.28]{Xu-Zhang}, one may see that
\[ Z(G) =\left\lbrace
  \begin{array}{ll}
    \lg a^{2^{r+s+v+t'-1}}\rg\lg b^{2^{s+t'}}\rg\cong C_{2}\times C_{2^{r+t-t'}} & s+t'\neq 0;\\
    \lg a^{2^{r+s+v+t'-1}}\rg\lg b^{2}\rg\cong C_{2}\times C_{2^{r+t-1}} & s+t'=0.
  \end{array} \right. \]
As $\lg a^{2^{r+s+v+t'-1}}\rg=\Omega_{1}(\lg a\rg)$, one has $Z(G)=\lg a^{2^{r+s+v+t'-1}}\rg$
due to $Z(G)$ is cyclic. It then from $tt'=0$ follows that $t'=r$ and $t=0$.
Therefore, $G$ has the following representation:
$$G=\lg a,b\ |\ a^{2^{2r+s+v}}=1, b^{2^{r+s}}=1, a^{b}=a^{-1+2^{r+v}}\rg,$$
where  $r\ge 2, sv=0.$\hfill\qed

\subsection{Two technical lemmas}

\begin{lem}\label{shu}
For any integer $n\ge 1$, if $0\le i\le n$,
then $2^{n-i}\mid \C_{2^{n}}^{i+1}.$
\end{lem}

\f\demo
Use induction on $n$. If $n=1$, then $i=0,1$, and in this case, a direct computation shows that $2^{n-i}=\C_{2^n}^{i+1}$.


Assume now that $n\ge 2$. If $i=n$, then $2^{n-i}=1$, and hence $2^{n-i}\mid \C_{2^{n}}^{i+1}.$
In what follows, assume that $i<n$, and that $2^{(n-1)-i}\mid \C_{2^{n-1}}^{i+1}.$
For any positive integer $m$ and for $0\leq i\leq 2^m-1$, let
\[\ell_m^i=(2^{m}-1)\cdots (2^{m}-(i+1)+1).\]
For any $1\le k\le m-1$, the highest power of $2$ that divides $(2^{m-1}-k)$ must also divide $(2^{m}-k)$.
Consequently, the highest power of $2$ that divides $\ell_{n-1}^i$ must also divide
$\ell_{n}^i$. As $\C_{2^{n-1}}^{i+1}=2^{n-1}\ell_{n-1}^i/(i+1)!$ and $\C_{2^{n}}^{i+1}=2^{n}\ell_{n}^{i}/(i+1)!$,
if the highest power of $2$ that divides $\C_{2^{n-1}}^{i+1}$ is $2^{t}$ then $2^{t+1}$ must also divide
$\C_{2^{n}}^{i+1}$. This completes the proof.\hfill\qed

\begin{lem}\label{2^k}
Let $G=\lg a,b\ |\ a^{2^{2r+s+t}}=1, b^{2^{r+s}}=1,a^{b}=a^{\pm 1+2^{r+t}}\rg$
with $r\ge 2, st=0.$ Let $1\neq g\in G$ be such that $\lg g\rg\cap\lg a\rg=1$ and $o(g)=2^{r+s-k}$ for some $0<k\leq r+s$.
Then there exists $b^{2^{k}}a^{j}\in G$ such that $\lg g\rg=\lg b^{2^{k}}a^{j}\rg$, where $2^{r+t+k}\mid j$.
\end{lem}

\f\demo For convenience of the statement, let $n=2^{2r+s+t}$, $m=2^{r+s}$ and $\g=\pm 1+2^{r+t}$. Then
$G=\lg a\rg: \lg b\rg\cong\mz_n: \mz_m$ and $b^{-1}ab=a^\g$. Let $g=b^{i'}a^{j'}$ for some $(i',j')\in\mz_m\times\mz_n$.
Since $\lg g\rg\cap\lg a\rg=1$, the pair $(i', j')$ is admissible. By Lemma~\ref{lem-admissible}, one has $o(b^{i'})=o(g)=2^{r+s-k}$, and hence
$2^{r+s-k}=\frac{2^{r+s}}{(2^{r+s}, i')}$. It follows that $(2^{r+s}, i')=2^k$, and hence there exists an odd integer $l$ such that
$$li'\equiv 2^{k}\ \ (\mod 2^{r+s}).$$
If $l>1$, then from Lemma~\ref{3} it follows that
$$g^{l}=(b^{i'}a^{j'})^{l}=b^{i'l}a^{j}=b^{2^{k}}a^{j},$$
for some $j\in\mz_{2^{2r+s+t}}.$
If $l=1$, then the above equation is also true.
Since $l$ is odd, one has
$$\lg g\rg=\lg g^{l}\rg=\lg b^{2^{k}}a^{j}\rg.$$

To complete the proof, it suffices to show that $2^{r+t+k}\ |\ j$. For convenience of the statement, let $i=2^k$ and $\ell=2^{r+s-k}$.
In view of $o(b^{2^{k}}a^{j})=o(g)=o(b^{2^k})=\ell$, by Lemma~\ref{lem-admissible-2}, one has
\begin{equation}\label{eq-lem-2^k-1}
j[\g^{i(\ell-1)}+\g^{i(\ell-2)}+\cdots+\g^i+1]\equiv 0\ (\mod n).
\end{equation}
If $\g^i=1$, then Eq.~(\ref{eq-lem-2^k-1}) implies that $j\ell\equiv 0\ (\mod n)$, and hence $2^{2r+s+t-(r+s-k)}=2^{r+t+k}$ divides $ j$.
Assume that $\g^i\neq 1$. Then Eq.~(\ref{eq-lem-2^k-1}) gives that \[j(\frac{\g^{2^{r+s}}-1}{\g^i-1})\equiv 0\ (\mod n).\]

Then to show that $2^{r+t+k}\ |\ j$, it suffices to prove the following claim.\medskip

\f{\bf Claim}\ {\em Let $d$ be a positive integer. Then $2^{r+t+d}$ is the highest power of $2$ that divides $ \g^{2^d}-1$.}

Recall that $\g=\pm 1+2^{r+t}$. Then \[\g^{2^d}-1=(\pm 1+2^{r+t})^{2^{d}}-1=\sum_{x=1}^{2^{d}}(\pm 1)^{x}\C_{2^{d}}^{x}\cdot 2^{(r+t)x}.\]

If $x>d+1$, then $2^{(r+t)x}>2^{(r+t)(d+1)}\geq 2^{r+t+d+1}$ since $r\ge 2$ and $t\ge 0$. Hence, $2^{r+t+d+1}\ |\ (\pm 1)^{x}\C_{2^{d}}^{x}\cdot 2^{(r+t)x}$ when $x>d+1$.  If $2\leq x\leq d+1$, then by Lemma~\ref{shu},  we have $2^{d-x+1}\mid \C_{2^{d}}^{x},$
and then
\[2^{d-x+1+(r+t)x}\ |\ (\pm 1)^{x}\C_{2^{d}}^{x}\cdot 2^{(r+t)x}.\]
Since $r\ge 2$ and $t\ge 0$, one has $2^{d-x+1+(r+t)x}\geq 2^{r+t+d+1}$. Again, we have $2^{r+t+d+1}\ |\ (\pm 1)^{x}\C_{2^{d}}^{x}\cdot 2^{(r+t)x}$ when $2\leq x\leq d+1$. If $x=1$, then $(\pm 1)^{x}\C_{2^{d}}^{x}\cdot 2^{(r+t)x}=\pm 2^{r+t+d}$. Thus, $2^{r+t+d}$ is the highest power of $2$ that divides $ \g^{2^d}-1$, completing the proof of our claim.
\hfill\qed

\subsection{Proof of Theorem \ref{Th}}

Let $G=\lg x\rg: \lg y\rg\cong \mz_{2^n}:\mz_{2^m}$.
If $G$ has a cyclic maximal subgroup, say $\lg a\rg$, then since $Z(G)$ is cyclic, one has $\Omega_{1}(Z(G))=\Omega_{1}(\lg a\rg)$.
Also, since $\lg x\rg\unlhd G$, one has $\Omega_{1}(Z(G))=\Omega_{1}(\lg x\rg)$.
If $o(y)>2$, then $1\neq y^{2}\in\lg a\rg$, and then $\Omega_{1}(Z(G))\leq\lg y\rg$, forcing $\lg x\rg\cap\lg y\rg\neq 1$, a contradiction.
Therefore, we have $o(y)=2$, and so $\lg x\rg$ is a cyclic maximal subgroup of $G$.
For any $1\neq g\in G$, if $\lg g\rg\cap\lg x\rg=1$, then $g^{2}\in\lg g\rg\cap\lg x\rg=1$. Hence $o(g)=2$. This implies that
$G$ is absolutely split.

In what follows, assume that $G$ has no cyclic maximal subgroups. By Lemma~\ref{G}, we may assume that
$$G=\lg a,b\ |\ a^{2^{2r+s+t}}=1, b^{2^{r+s}}=1,a^{b}=a^{\pm 1+2^{r+t}}\rg,$$
where $r\ge 2, st=0.$ Let $1\neq g\in G$ be such that $\lg g\rg\cap\lg a\rg=1$. Then $o(g)\leq o(b)$.
We may assume that $o(g)=2^{r+s-k}$ for some $1\leq k\leq r+s$.
By Lemma~\ref{2^k}, we have $\lg g\rg=\lg b^{2^k}a^{j}\rg$ with $2^{r+t+k}\mid j$. Without loss of generality,
we may assume that $g=b^{2^k}a^{j}$. Then $(2^k, j)\in\mz_{2^{2r+s+t}}\times\mz_{2^{r+s}}$ is a admissible pair with respect to $\lg a\rg$.

By the Claim in the proof of Lemma~\ref{2^k}, we have $2^{r+t+k}$ is the highest power of $2$ that divides $\g^{2^k}-1$, where $\g=\pm 1+2^{r+t}$.
Let $\frac{\g^{2^k}-1}{\g-1}=2^{u}\a$ for some odd integer $\a$. Then $2^u< 2^{r+t+k}$ and then there exists $\a^{-1}\in\mz_{2^{2r+s+t}}$ such that $\a^{-1}\cdot\a\equiv1\ (\mod 2^{2r+s+t})$. Let $l\in\mz_{2^{2r+s+t}}$ be such that $l\equiv \a^{-1}\frac{j}{2^{u}}\ (\mod 2^{2r+s+t})$.
Then \[l(\frac{\g^{2^k}-1}{\g-1})=l\cdot 2^{u}\cdot \a\equiv j\ (\mod 2^{2r+s+t}).\]
Now to complete the proof, by Theorem~\ref{thm-absolute-1}, it suffices to show that the pair $(1,l)$ is admissible.

Actually, we have
\[
\begin{array}{lll}
l[\g^{(2^{r+s}-1)}+\g^{(2^{r+s}-2)}+\cdots+\g+1]&=&l(\frac{\g^{2^{r+s}}-1}{\g-1})\\
&=&\frac{\g^{2^{r+s}}-1}{\g^{2^k}-1}\cdot l(\frac{\g^{2^{k}}-1}{\g-1})\\
&\equiv& \frac{\g^{2^{r+s}}-1}{\g^{2^k}-1}\cdot j\ (\mod 2^{2r+s+t}).\\
\end{array}\]
By the Claim in the proof of Lemma~\ref{2^k}, we have $2^{r+s-k}\ |\ \frac{\g^{2^{r+s}}-1}{\g^{2^k}-1}$, and since $2^{r+t+k}\mid j$, one has
$\frac{\g^{2^{r+s}}-1}{\g^{2^k}-1}\cdot j\equiv 0\  (\mod 2^{2r+s+t}).$
Thus, \[l[\g^{(2^{r+s}-1)}+\g^{(2^{r+s}-2)}+\cdots+\g+1]\equiv 0\  (\mod 2^{2r+s+t}).\]
By Lemma~\ref{lem-admissible-2}, the pair $(1,l)$ is admissible.
\hfill\qed

\section{Metacirculants of $2$-power order}

\begin{lem}\label{O}
Let $G$ be a split metacyclic $2$-group. If the center of $G$ is not cyclic,
then $\Omega_{1}(G)\cong C_2\times C_2.$
\end{lem}

\f\demo Let $G=\lg a\rg:\lg b\rg\cong C_{2^{n}}: C_{2^{m}}$.
Then $$C_{2}\times C_{2}\cong \lg a^{2^{n-1}}\rg\times\lg b^{2^{m-1}}\rg\leq \Omega_{1}(G).$$
Since $Z(G)$ is noncyclic, one has $\lg a^{2^{n-1}}\rg\times\lg b^{2^{m-1}}\rg\leq Z(G)$.
For every $g\in \Omega_{1}(G)$, let $g=b^{i}a^{j}$ for some $(i,j)\in\mz_{2^m}\times\mz_{2^n}$.
Then $g^{2}=1$. Applying Lemma~\ref{3}, we obtain that $b^{2i}=1$, and then we have
$b^{i}\in \lg b^{2^{m-1}}\rg\leq Z(G)$. It then follows that $1=g^{2}=(b^{i}a^{j})^{2}=b^{2i}a^{2j}=a^{2j}.$
Therefore $a^{j}\in \lg a^{2^{n-1}}\rg\leq Z(G)$.
Therefore $g\in  \lg a^{2^{n-1}}\rg\times\lg b^{2^{m-1}}\rg$, and hence $\Omega_{1}(G)= \lg a^{2^{n-1}}\rg\times\lg b^{2^{m-1}}\rg\cong C_2\times C_2$.
\hfill\qqed

The following theorem shows that the answer to Question~A is positive in the case when the graph under question has order a $2$-power.

\begin{thm}\label{Result}
A connected weak metacirculant with order a $2$-power is a metacirculant if and only if it is a split weak metacirculant.
\end{thm}

\f\demo By \cite[Lemma~2.2]{LWS}, it suffices to prove the sufficiency. Let $\G$ be a split weak metacircualnts of order a $2$-power. By the definition of split weak metacirculant, ${\rm Aut(\Gamma)}$ has a split metacyclic subgroup $X$ which is transitive on $V(\G)$. Let $G$ be a Sylow $2$-subgroup of $X$. Then $G$ is also split metacyclic, and by \cite[Theorem 3.4]{20},
$G$ is also transitive on $V(\G)$.
If $G$ is regular on $V(\Gamma)$, then $\Gamma$ is a Cayley graph on $G$, and then $\Gamma$ must be a metacirculant, as required.
In what follows, we always assume that $G$ is not regular on $V(\Gamma)$.

If $Z(G)$ is non-cyclic, then by Lemma~\ref{O}, we have $\Omega_{1}(G)\cong C_{2}\times C_{2}$ and $\Omega_{1}(G)\leq Z(G)$.
For any $v\in V(\Gamma)$, if $G_{v}\neq 1$, then $G_{v}\cap \Omega_{1}(G)\neq 1$.
However, $G_{v}\cap \Omega_{1}(G)\unlhd G$, forcing that $G_{v}\cap \Omega_{1}(G)$ fixes every vertex of $\Gamma$, a contradiction.
Thus, $G_{v}=1$, and so $G$ is regular on $V(\Gamma)$.
It follows that $\Gamma$ is a Cayley graph on $G$.  Again, this is impossible.

Assume now $Z(G)$ is cyclic. By Theorem~\ref{Th}, $G$ is absolutely split. We may assume that $G=\lg x\rg:\lg y\rg\cong C_{2^{n}}:C_{2^m}$, and that $G$ is absolutely split with respect to $\lg x\rg$. Since $G$ is transitive on $V(\Gamma)$, $\lg x\rg$ acts semiregularly on $V(\Gamma)$.
Assume that $\lg x \rg$ has $2^{\ell}$ orbits for some $1\leq\ell\leq n$. Then the kernel of $G$ acting on the orbits of $\lg x\rg$ is $\lg x\rg: G_{v}$ for some $v\in V(\G)$. So $y^{2^{\ell}}\in \lg x\rg: G_{v}$ and $|G_{v}|=2^{m-\ell}$.
Let $G_{v}=\lg z\rg$.
Since $\lg x \rg\cap G_{v}=1$,
one has $\lg z\rg\cap\lg x\rg=1$.
Since $G$ is absolutely split, there exists $y'\in G$ such that $z\in\lg y'\rg\cong C_{2^{m}}$.
Moreover, $G=\lg x\rg: \lg y'\rg$.
Then $y'$ cyclically permutes the $2^{\ell}$ orbits of $\lg x\rg$,
and $(y')^{2^{\ell}}\in G_{v}$.
This implies that $\Gamma$ is a metacirculant.
\hfill\qed

\medskip
\f {\bf Acknowledgements:}\ This work was partially supported by the National
Natural Science Foundation of China (11271012, 11671030) and the Fundamental
Research Funds for the Central Universities (2015JBM110).

\end{document}